\documentclass{article}
\usepackage{graphicx}
\usepackage{amsmath}

\input{tcilatex}

\begin{document}

\title{On fixed points of Poisson shot noise transforms}
\author{Aleksander M. Iksanov and Zbigniew J. Jurek \\
Cybernetics Faculty, Kiev T.Shevchenko, \\
National University, 01033 Kiev, Ukraine \\
Department of Mathematics, Wayne State University\\
Detroit, MI 48202, USA\\
iksan@unicyb.kiev.ua and zjjurek@math.wayne.edu}
\date{\textbf{Published} in Advances in Applied Probability, 34, no.4, pp.798-825.}
\maketitle

\begin{abstract}
Distributional fixed points of a Poisson shot noise transform (for
nonnegative and nonincreasing response functions bounded by 1) are
characterized. The tail behavior of fixed points is described. Typically
they have either exponential moments or their tails are proportional to a
power function, with exponent greater than $-1$. The uniqueness of fixed
points is also discussed. Finally, it is proved that in most cases fixed
points are absolutely continuous, apart from the possible atom at zero.

\end{abstract}

\end{document}